\documentclass[10pt]{article}

\usepackage{a4wide}
\usepackage{amssymb}
\usepackage{amsfonts}
\usepackage{amsmath}
\input xy
\xyoption{arrow} \xyoption{matrix}

\date{}

\newtheorem{proposition}{Proposition}[section]
\newtheorem{theorem}[proposition]{Theorem}
\newtheorem{lemma}[proposition]{Lemma}

\newtheorem{corollary}[proposition]{Corollary}

\def\der{\partial }

\def\nFM0{{\nu }_{F,M_0}}
\def\nFN0{{\nu }_{F,N_0}}
\def\nGN0{{\nu }_{G,N_0}}

\def\N0{ {\bf N}_0 }

\def\t{\otimes}
\def\g{\gamma}

\def\ra{\rightarrow}

\def\Xpm{X^{\pm }}

\def\s{\sigma}
\def\Z{\mathbb{Z}}

\def\l1{{\lambda}_1}

\def\a{\alpha}
\def\a0{ {\alpha }_0}
\def\a1{ {\alpha }_1}

\def\l{\lambda}

\def\nFGM0{{\nu }_{F,G,M_0}}

\def\nFN0{{\nu}_{F,N_0}}

\def\sm{{\sigma}^m}

\def\sm1{{\sigma}^{-1}}

\def\smtp1{{\sigma}^{-t+1}}

\def\S1{S^{-1}}

\def\Xpm1{X^{\pm 1}_1}

\def\sPM1{{\sigma }^{\pm 1}}
\def\sMP1{{\sigma }^{\mp 1 }}

\def\d{\delta}

\def\di{{\rm d.ind}}

\def\L{\Lambda}

\def\Ytm1{Y^{t-1}}
\def\Yim1{Y^{i-1}}

\def\CK{{\cal K}}

\def\Aut{{\rm Aut}}

\def\bA{\overline{A}}

\def\ad{{\rm ad }}

\def\char{{\rm char }}

\def\ker{ {\rm ker } }

\def\SL2Z{ {\rm SL}_2({\bf Z}) }

\def\Gp1{ G^{1 , 1 } }
\def\P11{ P^{-1 , 1 } }
\def\Pp1{ P^{1 , 1 } }

\def\nCLsr{{}^\nu\kern-2pt {\cal L}^{\sigma , \rho  }}
\def\nP{{}^\nu \kern-2pt P}
\def\nL{{}^\nu\kern-2pt L}
\def\nLL{{}^\nu\kern-2pt \Lambda}
\def\nPsr{{}^\nu\kern-2pt P^{\sigma , \rho  }}
\def\nLsr{{}^\nu\kern-2pt L^{\sigma , \rho  }}
\def\nuCL{{}^\nu\kern-2pt  {\cal L}}
\def\nCLsr{{}^\nu\kern-2pt {\cal L}^{\sigma , \rho  }}
\def\nCL1m{{}^\nu\kern-2pt {\cal L}^{-1 , 1  }}
\def\x1nu{x^\frac{1}{\nu}}
\def\xm1nu{x^{-\frac{1}{\nu}}}

\def\ra{\rightarrow }

\def\nAM0{{\nu }_{{\cal A},M_0}}
\def\nAN0{{\nu }_{{\cal A},N_0}}

\def\ad{ {\rm ad }}

\def\ga{\mathfrak{a}}

\def\SL{{\rm SL}}

\def\di!{\frac{\der^i}{i!}}
\def\dik!{\frac{\der^k_i}{k!}}
\def\N{\mathbb{N}}

\def\0{\overline{0}}
\def\1{\overline{1}}

\def\Ln1{\L_{n,\overline{1}}}

\def\a1{a_{\overline{1}}}

\def\S{\Sigma}

\def\vn1{\overrightarrow{n-1}}

\def\Q{\mathbb{Q}}

\def\mS{\mathbb{S}}

\def\clKdim{{\rm cl.Kdim}}
\def\pd{{\rm pd}}

\def\lgldim{{\rm l.gldim}}
\def\rgldim{{\rm r.gldim}}

\def\K1{{\rm K}_1}

\def\wdim{{\rm w.dim}}
\def\btaut{\overline{\tau_t}}

\begin{document}

\author{V. V. \  Bavula %(HilbSn.tex)
}

\title{An analogue of  Hilbert's Syzygy Theorem for the algebra
of one-sided inverses of a polynomial algebra}

\maketitle

\begin{abstract}
An analogue of  Hilbert's Syzygy Theorem is proved for the algebra
$\mS_n (A)$ of one-sided inverses of the  polynomial algebra
$A[x_1, \ldots , x_n]$ over an arbitrary ring $A$:
$$ \lgldim (\mS_n(A))= \lgldim (A) +n.$$
The algebra $\mS_n(A)$ is noncommutative,  neither left nor right
Noetherian and not a domain.  The proof is based on a
generalization of the Theorem of Kaplansky (on the projective
dimension) obtained in the paper. As a consequence it is proved
that  for a left or right Noetherian algebra $A$:
$$ \wdim (\mS_n(A))= \wdim (A) +n.$$

%$\noindent $



%$${\bf Contents}$$
%\begin{enumerate}
%\item Introduction. \item  \item The group of automorphisms of the
%algebra $\mS_2$.
%\end{enumerate}
\end{abstract}

%%%%%%%%%%%%%%%%%% SECTION 1 %%%%%%%%%%%%%%%%%%%%%%%%

\section{Introduction}
 Throughout, ring
means an associative ring with $1$; module means a left module;
 $\N :=\{0, 1, \ldots \}$ is the set of natural numbers.
%; $K$ is a field  and  $K^*$ is its group of units; $P_n:= K[x_1,
%\ldots , x_n]$ is a polynomial algebra over $K$;
%$\der_1:=\frac{\der}{\der x_1}, \ldots , \der_n:=\frac{\der}{\der
%x_n}$ are the partial derivatives ($K$-linear derivations) of
%$P_n$; $\End_K(P_n)$ is the algebra of all $K$-linear maps from
%$P_n$ to $P_n$ and $\Aut_K(P_n)$ is its group of units (i.e. the
%group of all the invertible linear maps from $P_n$ to $P_n$); the
%subalgebra $A_n:= K \langle x_1, \ldots , x_n , \der_1, \ldots ,
%\der_n\rangle$ of $\End_K(P_n)$ is called the $n$'th {\em Weyl}
%algebra.

$\noindent $

{\it Definition}, \cite{shrekalg}. Let $P_n(A)=A[x_1, \ldots ,
x_n]$ be a polynomial ring with coefficients in a ring $A$.
 The %{\em Shrek}
{\em algebra} $\mathbb{S}_n=\mS_n(A)$ {\em of one-sided inverses}
of $P_n(A)$ is a ring  generated over  the ring  $A$ of
 by $2n$ elements $x_1, \ldots , x_n, y_n,
\ldots , y_n$ that satisfy the defining relations:
$$ y_1x_1=1, \ldots , y_nx_n=1, \;\; [x_i, y_j]=[x_i, x_j]= [y_i,y_j]=0
\;\; {\rm for\; all}\; i\neq j,$$ where $[a,b]:= ab-ba$ is  the
commutator of elements $a$ and $b$.

$\noindent $

By the very definition, the ring $\mS_n$ is obtained from the
polynomial ring $P_n$ by adding commuting, left (but not
two-sided) inverses of its canonical generators. When $A=K$ is a
field, the algebra $\mS_1=K\langle x, y\, | \, yx=1\rangle$ is a
well-known primitive algebra \cite{Jacobson-StrRing}, p. 35,
Example 2.  Over the field
 $\mathbb{C}$ of complex numbers, the completion of the algebra
 $\mS_1$ is the {\em Toeplitz algebra} which is the
 $C^*$-algebra generated by a unilateral shift on the
 Hilbert space $l^2(\N )$ (note that $y=x^*$). The Toeplitz
 algebra is the universal $C^*$-algebra generated by a
 proper isometry. If $\char (K)=0$ then the algebra $\mS_n$ is
 isomorphic to the algebra $K\langle
\frac{\der}{\der x_1}, \ldots ,\frac{\der}{\der x_n},  \int_1,
\ldots , \int_n\rangle $ of  scalar  integro-differential
operators (via $x_i \mapsto \int_i$, $y_i\mapsto \frac{\der}{\der
x_i}$).

$\noindent $

The algebra $\mS_n=\mS_n(K)$ was studied in detail in
\cite{shrekalg}: the Gelfand-Kirillov dimension of the algebra
$\mS_n$ is $2n$, $\clKdim (\mS_n)=2n$,  the weak and the global
dimensions of $\mS_n$ are $n$. The algebra $\mS_n$ is neither left
nor right Noetherian as was shown by Jacobson
\cite{Jacobson-onesidedinv-1950} when $n=1$ (see also Baer
\cite{Baer-inverses-1942}). Moreover, it contains infinite direct
sums of left and right ideals. In the series of three papers
\cite{shrekaut}, \cite{K1aut} and \cite{Snaut} the group $G_n:=
\Aut_{K-{\rm alg}}(\mS_n)$  of automorphisms and the group
$\mS_n^*$ of units of the algebra $\mS_n$ and there explicit
generators were found (both groups are huge).  The group $G_1$ was
found by Gerritzen, \cite{Gerritzen-yx=1-2000}). In \cite{K1aut}
and \cite{K1-Sn-group} it is proved that ${\rm K}_1(\mS_1) \simeq
K^*$ and  ${\rm K}_1(\mS_n) \simeq K^*$ for all $n>1$
respectively.

 The aim of the
present paper is to prove an analogue of  Hilbert's Syzygy Theorem
for the ring $\mS_n(A)$ over an arbitrary ring $A$.

$\noindent $

{\em   Hilbert's Syzygy Theorem} \cite{ERZ},
\cite{Bass-book-K-theory} states that for any ring $A$,
$$ \lgldim (A[x_1, \ldots , x_n])= \lgldim (A) +n.$$
The reason  why it holds is the following Theorem of Kaplansky.
\begin{theorem}\label{KaplanskyTHM}%\marginpar{KaplanskyTHM}
{\rm \cite{Kaplansky-book-69}} Let $A$ be a ring, $s$ be its
regular and central element, $\bA := A/(s)$. If $M$ is a nonzero
$\bA$-module with $\pd_{\bA}(M)=n<\infty$ then $\pd_{A}(M)=n+1$.
\end{theorem}

The aim of the paper is to prove an analogue of Hilbert's Syzygy
Theorem for the ring $\mS_n(A)$.

\begin{theorem}\label{AH22May10}%\marginpar{AH22May10}
 Let $A$ be a ring. Then
$$ \lgldim (\mS_n(A))= \lgldim (A) +n \;\; {\rm and}\;\;  \rgldim (\mS_n(A))= \rgldim (A) +n.$$
\end{theorem}

The proof of this theorem is based on the following generalization
of the Theorem of Kaplansky proved in the paper.

\begin{theorem}\label{22May10}%\marginpar{22May10}
 Let $A$ be a ring, $s$ be a regular and left normal
element of the ring $A$, $t$ be a right regular and right normal
element of the ring $A$ such that $(s)=(t)$, and  $\bA := A/(s)$.
Suppose that  either
\begin{enumerate}
\item  $ta-at\in (t^2)$ for all elements $a\in A$, or \item the
ring endomorphism $\btaut$ of the ring $\bA$ is an automorphism
(see (\ref{ttt1})).
\end{enumerate}
If $M$ is a nonzero $\bA$-module with $\pd_{\bA}(M)=n<\infty$ then
$\pd_{A}(M)=n+1$.
\end{theorem}

The Theorem of Kaplansky follows from Theorem \ref{22May10} when
we set $t=s$ (both conditions hold).

$\noindent $

Analogues of Hilbert's Syzygy Theorems for some (mostly
Noetherian) rings under additional assumptions on the coefficient
ring $A$ (like being Noetherian) were given in \cite{Aus},
\cite{ERZ}, \cite{Bav-gldimGWA-1994} and \cite{Bav-THM-96}.
 Generalizations of the Theorem of Kaplansky are the key in
 finding the global dimension of certain classes of Noetherian
 rings (like skew polynomial rings, skew Laurent polynomial rings,
 rings of differential operators, etc) see \cite{Bjork-72},
 \cite{Goodearl-75}, \cite{Rin-Ros-76}, \cite{Ros-Stafford-76} and
 \cite{Brown-Haj-M-82}.

%%%%%%%%%%%%%%%%%% SECTION 2  %%%%%%%%%%%%%%%%%%%%%%%%

\section{Proofs of Theorems \ref{AH22May10} and \ref{22May10}}\label{}%\marginpar{}

At the beginning of the  section we collect some basic facts on
the ring $\mS_n(A)$ and prove several lemmas that are used in the
proofs of Theorems \ref{AH22May10} and \ref{22May10}.

$\noindent $

{\bf The ring $\mS_n(A)$  of one-sided inverses of a polynomial
ring}. Let $A$ be a ring.  The ring $\mS_n(A)$ is the direct sum
of $\N^2$ copies of ${}_AA$ and $A_A$, %\marginpar{SnAS}
\begin{equation}\label{SnAS}
\mS_n(A)=A\t_\Z \mS_n(\Z )=\bigoplus_{\alpha , \beta \in \N^n}
Ax^\alpha y^\beta =\bigoplus_{\alpha , \beta \in \N^n} x^\alpha
y^\beta A
\end{equation}
where $x^\alpha := x_1^{\alpha_1} \cdots x_n^{\alpha_n}$, $\alpha
= (\alpha_1, \ldots , \alpha_n)$, $y^\beta := y_1^{\beta_1} \cdots
y_n^{\beta_n}$ and $\beta = (\beta_1,\ldots , \beta_n)$.
 The ring $\mS_n(A)$ contains the
ideal %\marginpar{1SnAs}
\begin{equation}\label{1SnAs}
F(A):=\bigoplus_{\alpha , \beta \in \N^n}AE_{\alpha \beta}, \;\;
{\rm where}\;\; E_{\alpha \beta}:=\prod_{i=1}^n E_{\alpha_i
\beta_i}(i), \;\; E_{\alpha_i
\beta_i}(i):=x_i^{\alpha_i}y_i^{\beta_i}-x_i^{\alpha_i+1}y_i^{\beta_i+1}.
\end{equation}
Note that $E_{\alpha \beta}E_{\g \rho}=\d_{\beta \g }E_{\alpha
\rho}$ for all elements $\alpha, \beta , \g , \rho \in \N^n$ where
$\d_{\beta
 \g }$ is the Kronecker delta function. Let $\ga_n(A)$ be the ideal
 of the algebra $\mS_n (A)$ generated by the elements $E_{00}(i)$,
 $i=1, \ldots , n$ (if $n=1$ then $\ga_1(A) = F(A)$).
 The factor ring $\mS_n(A)
 / \ga_n(A)$ is canonically isomorphic to the Laurent polynomial
 $A$-algebra $L_n (A) := A[ x_1, x_1^{-1}, \ldots , x_n ,
 x_n^{-1}]$ via the $A$-isomorphism:
 %\marginpar{LnSn}
\begin{equation}\label{LnSn}
\mS_n(A) / \ga_n(A)\ra L_n(A), \;\; x_i\mapsto x_i, \;\;
y_i\mapsto x_i^{-1}, \;\; i=1, \ldots , n,
\end{equation}
since $y_ix_i=1$ and $1-x_iy_i= E_{00}(i) \in \ga_n(A)$.

Suppose that the ring $A$ admits an involution $*$, i.e. it is an
anti-isomorphism ($(ab)^*= b^*a^*$ for all elements $a,b\in A$)
with $a^{**}=a$ for all elements $a\in A$. The involution $*$ can
be extended to the involution $*$ on the ring $\mS_n(A)$ by the
rule
$$ x_i^* = y_i, \;\; y_i^*= x_i, \;\; i=1, \ldots , n.$$
In particular when  $A$ is a commutative ring and $*$ is the
identity map on $A$ we have the $A$-involution (i.e. $a^*=a$ for
all $a\in A$) on the ring $\mS_n(A)$. There is a split ring
$A$-epimorphism %\marginpar{pisplit}
\begin{equation}\label{pisplit}
\pi : \mS_n(A) \ra A, \;\; x_i\mapsto 1, \;\; y_i\mapsto 1, \;\;
i=1, \ldots , n.
\end{equation}
Therefore, %\marginpar{SA1}
\begin{equation}\label{SA1}
\mS_n(A)=A\bigoplus (x_1-1, \ldots , x_n-1, y_1-1, \ldots ,
y_n-1).
\end{equation}
Lemma \ref{a22May10}.(2) shows that the ideal $(x_1-1, \ldots ,
x_n-1, y_1-1, \ldots , y_n-1)$ of the ring $\mS_n(A)$  is equal to
the ideal $(x_1-1, \ldots , x_n-1)$ or $(y_1-1, \ldots , y_n-1)$
or $(t_1-1, \ldots , t_n-1)$ where $t_i = x_i, y_i$.

 The polynomial ring  $P_n(A) = A[x_1, \ldots , x_n]$
 is an $\mS_n(A)$-module. In more detail,
 $${}_{\mS_n(A)}P_n(A)\simeq \mS_n(A) / (\sum_{i=0}^n \mS_n(A)
y_i) =\bigoplus_{\alpha \in \N^n} Ax^\alpha \overline{1},\;\;
\overline{1}:= 1+\sum_{i=1}^n \mS_n(A)y_i,$$ and the action of the
canonical generators of the $A$-algebra $\mS_n(A)$ on the
polynomial ring $P_n(A) $ is given by the rule:
$$ x_i*x^\alpha = x^{\alpha + e_i}, \;\; y_i*x^\alpha = \begin{cases}
x^{\alpha - e_i}& \text{if } \; \alpha_i>0,\\
0& \text{if }\; \alpha_i=0,\\
\end{cases}  \;\; {\rm and }\;\; E_{\beta \g}*x^\alpha = \d_{\g
\alpha} x^\beta,
$$
where the set $e_1:= (1,0,\ldots , 0),  \ldots , e_n:=(0, \ldots ,
0,1)$ is the canonical basis for the free $\Z$-module
$\Z^n=\bigoplus_{i=1}^n \Z e_i$.

An element $r$ of a ring $R$ is called a {\em left normal} if $Rr=
(r)$, i.e. $rR\subseteq Rr$, where $(r)=RrR$ is the (two-sided)
ideal of the ring $R$ generated by the element $r$. An element $r$
of a ring $R$ is called a {\em right normal} if $rR= (r)$, i.e.
$Rr\subseteq rR$. An element $r$ of a ring $R$ is called a {\em
normal} if it is left and right normal, equivalently, $Rr=rR$
(since $Rr= (r) = rR$).

\begin{lemma}\label{a22May10}%\marginpar{a22May10}
Let $A$ be a ring. Then
\begin{enumerate}
\item The element $x-1\in \mS_1(A)$ is left normal and the element
$y-1\in \mS_1(A)$ is right normal, i.e. $\mS_1(A) (x-1)= (x-1)$
and $(y-1)\mS_1(A) = (y-1)$.  \item  $(x-1) = (y-1)$,
${}_A\mS_1(A)_A= A\bigoplus (x-1)$, and $F(A) \subseteq (x-1)$.
\end{enumerate}
\end{lemma}

{\em Remark}. Lemma \ref{b22May10} shows that the element $x-1$ is
not a right normal element and the element $y-1$ is not a left
regular of $\mS_n(A)$.

$\noindent $

{\it Proof}. By (\ref{SnAS})  and since  $x-1\in \mS_1(\Z )$ and
$y-1\in \mS_1(\Z )$, it suffices to prove the lemma for $A=\Z$.
The $\Z$-algebra $\mS_1=\mS_1(\Z )$ is generated by the elements
$x$ and $y$ over the ring $\Z$. The left ideal $\mS_1(x-1)$ is an
ideal of the ring $\mS_1$ since %\marginpar{x1y}
\begin{equation}\label{x1y}
(x-1)y = (1-(x-1)y)(x-1).
\end{equation}
In more detail, using the equalities $E_{00} = 1-xy$ and
$E_{00}x=0$, we see that $$(1-(x-1)y)(x-1)= (E_{00}+y)(x-1) =
-E_{00}+1-y= xy-1+1-y= (x-1)y.$$ By (\ref{x1y}), the left ideal
$\mS_1(x-1)$ is an ideal of $\mS_1$, i.e. $x-1$ is a left normal
element of the ring $\mS_1$. The ring $\mS_1$ admits the
involution $*$ over the ring $\Z$: $y^*=x$. Then the element
$y-1=(x-1)^*$ is right normal. Notice that %\marginpar{x2y}
\begin{equation}\label{x2y}
1-y= y(x-1)\in (x-1).
\end{equation}
Then $(y-1) \subseteq (x-1)$ and the opposite inclusion follows
when we apply the involution $*$ to the original one, i.e. $(x-1)=
(y-1)$. Now, %\marginpar{S1Z}
\begin{equation}\label{S1Z}
\mS_1 = \Z [x]\t_\Z \Z [y]= (\Z \bigoplus \Z [x] (x-1))\t_\Z  (\Z
\bigoplus \Z [y] (y-1))= \Z \bigoplus (x-1)
\end{equation}
 since the
$\Z$-algebra epimorphism $\pi : \mS_1\ra \Z$, $ x\mapsto 1$, $
y\mapsto 1$, is a split epimorphism. Since $\ker (\pi )= (x-1)$
and $\pi (E_{ij})=0$ for all elements $i,j\in  \N$, we have the
inclusion $F\subseteq (x-1)$.  $\Box $

%$\noindent $

\begin{lemma}\label{b22May10}%\marginpar{b22May10}
Let $A$ be a ring. Then neither the right ideal $(x-1)\mS_1(A)$
nor the left ideal $\mS_1(A)(y-1)$ is an  ideal of the ring
$\mS_1(A)$.
\end{lemma}

{\it Proof}. By (\ref{SnAS})  and since  $x-1, y-1\in \mS_1(\Z)$,
we may assume that $A=\Z$. In view of the involution $*$ on the
ring  $\mS_1= \mS_1(\Z )$ it suffices to show that the right ideal
$(x-1)\mS_1$ is not an ideal of the ring $\mS_1$ since
$\mS_1(y-1)=((x-1)\mS_1)^*$. Suppose that $(x-1)\mS_1$ is an ideal
of the ring $\mS_1$, we seek a contradiction. Since $F(\Q ) =\Q
\t_\Z F$ (where $F=F(\Z )$) is the least non-zero ideal of the
algebra $\mS_1(\Q )$ (Proposition 2.5.(2), \cite{shrekalg}) and
$\mS_1\subseteq \mS_1(\Q )$, it follows from  (\ref{1SnAs}) that
$F(\Q) \cap (x-1)\mS_1= m F$ for a nonzero integer $m$ since $F(\Q
) \subseteq (x-1)\mS_1(\Q )$.   Since $mE_{00}\in F(\Q )\cap
(x-1)\mS_1$ there exists an element $a\in \mS_1$ such that
$(x-1)a=mE_{00}$. Taking this equality modulo $F$ and using the
fact that the factor ring $\mS_1/F$ is a Laurent polynomial
algebra $\Z [ x,x^{-1}]$ over $\Z$ which is a domain, we see that
$a\in F$. The ideal $F$ is the direct sum $\bigoplus_{j\in
\N}E_{\N , j}$ of left ideals $E_{\N , j}=\bigoplus_{i\in \N}\Z
E_{ij}\simeq {}_{\mS_1}\Z[x]$ ($E_{0j}\mapsto 1$). Recall that
$E_{ij}=x^iy^j-x^{i+1}y^{j+1}$. Since the left multiplication by
the element $x-1$ in $\Z [x]$ is an injection and $E_{00}\in E_{\N
, 0}$, we must have $a\in E_{\N , 0}$. Notice that ${}_{\Z
[x]}E_{\N, 0}\simeq {}_{\Z [x]}\Z [x]$, $E_{00}\mapsto 1$. The
equality $(x-1)a=mE_{00}$ implies that $m\in (x-1)\Z [x]$, a
contradiction. $\Box $

$\noindent $

An element $r$ of a ring $R$ is called a {\em left  regular} if
$sr=0$ implies $s=0$ for $s\in R$. Similarly, {\em right regular}
is defined, and {\em regular} means both left  and right regular.

\begin{lemma}\label{c22May10}%\marginpar{c22May10}
Let $A$ be a ring. The elements $x-1$ and $y-1$ of the algebra
$\mS_1(A)$ are regular.
\end{lemma}

{\it Proof}. By (\ref{SnAS})  and since $x-1, y-1\in \mS_1(\Z )$,
we may assume that $A=\Z$. In view of the involution $*$ on the
algebra $\mS_1=\mS_1(\Z )$ ($y^*=x$), it suffices to show that the
element $x-1$ is regular. Suppose that $l(x-1)=0$ and $(x-1)r=0$
for some elements $l,r\in \mS_1$. We have to show that $l=r=0$.
Since the element $x-1$ is a regular element of the factor ring
$\mS_1/F = \Z [x,x^{-1}]$, we see that $l,r\in F$. We have seen in
the proof of Lemma \ref{b22May10} that ${}_{\Z [x]}F\simeq \Z
[x]^{(\N )}$ hence $r=0$. The ideal $F$ is the direct sum
$\bigoplus_{i\in \N}E_{i,\N}$ of right ideals
$E_{i,\N}=\bigoplus_{j\in \N}\Z E_{ij}\simeq \mS_1/x\mS_1\simeq \Z
[y]_{\mS_1}$ ($E_{i0}\mapsto 1$). In the right $\mS_1$-module $\Z
[y]$, $y^ix=\begin{cases}
y^{i-1}& \text{if }i>0,\\
0& \text{if }i=0, \\
\end{cases}$
 i.e. the map $\cdot x: \Z [x]\ra \Z [x]$, $v\mapsto vx$, is a
 locally nilpotent map, that is $\Z [y] =\bigcup_{i\geq 1}\ker_{\Z
 [y]}(\cdot x^i)$. Therefore, the map $x-1$ an isomorphism, hence
 $l=0$. This proves that the element $x-1$ is regular. $\Box $

$\noindent $

Let $A$ be a ring, $\s$ be its ring endomorphism, and $M$ be an
$A$-module. The {\em twisted} by $\s$ $A$-module ${}^\s M$ is
equal to $M$ as an abelian group but the $A$-module structure on
${}^\s M $ is given by the rule $a\cdot m = \s (a) m$ where $a\in
A$ and $m\in M$.

$\noindent $

Any right regular and right normal element $t$ of a ring $A$
determines the ring monomorphism $\tau_t : A\ra A$ by the rule
%\marginpar{ttt}
\begin{equation}\label{ttt}
at=t\tau_t(a)\;\; {\rm for }\; \; a\in A.
\end{equation}
Clearly, $\tau_t(t) = t$, and so $\tau_t((t))=(t)$ and we have the
ring endomorphism $\btaut : \bA:=A/(t)\ra \bA$ given by the rule
%\marginpar{ttt1}
\begin{equation}\label{ttt1}
\btaut (a+(t))=\tau_t(a) +(t).
\end{equation}

$\noindent $

{\bf  Proof of Theorem \ref{22May10}}. Since the element $s$ is
left regular ($as=0$ implies $a=0$ for $a\in A$) and left normal
($As= (s)$) there is the short exact sequence of $A$-modules
$$
0\ra A \stackrel{\cdot s}{\ra} A\ra \bA \ra 0.$$  This implies
that $\pd_A(P)\leq 1$  for each projective $\bA$-module $P$. In
fact, $\pd_A(P)=1$ as the $A$-module $P$ cannot be projective
since $sP=0$ and the element $s\in A$ is right regular. This
proves the theorem when $n=0$. Let $n>0$ and we use induction on
$n$. Let $0\ra N \ra P\ra M\ra 0$ be a short exact sequence of
$\bA$-modules where $P$ is a projective $\bA$-module. Then $N\neq
0$ and $\pd_{\bA}(N)=n-1$. By induction $\pd_A(N) =n$ and
$\pd_A(P)=1$. Therefore, $\pd_A(M)\leq n+1$ and $\pd_A(M) = n+1$
provided $n>1$.

If $n=1$ then $M = Q/R$ for some {\em free} $A$-module $Q$ where
$R$ is its $A$-submodule such that $(s)Q\subseteq R$. Then the
following short sequences of $\bA$-modules  are exact:
%\marginpar{RQM}
\begin{equation}\label{RQM}
0\ra R/(s)Q\ra Q/(s)Q\ra M\ra 0,
\end{equation}
%\marginpar{1RQM}
\begin{equation}\label{1RQM}
0\ra (s)Q/(s)R\ra R/(s)R\ra R/(s)Q\ra 0.
\end{equation}
By the assumption, $(s) = (t) = tA$ and $At\subseteq tA$. Then
$(s)Q/(s)R=tQ/tR$. Since the element $t\in A$ is right regular,
the map %\marginpar{2RQM}
\begin{equation}\label{2RQM}
tQ/tR\ra Q/R, \;\; tq+tR\mapsto q+R,
\end{equation}
is a bijection. We claim that $$\pd_{\bA}((s)Q/(s)R)=1.$$ Suppose
that the first condition in the theorem holds. Since $at-ta\in
(t^2)$ for all elements $a\in A$ and $(t)R=tR$, the map
(\ref{2RQM}) is an $A$-module (and $\bA$-module) isomorphism:
$$ a(tq+(t)R)=atq+(t)R=taq+ (at-ta)q+(t)R=taq +(t)R$$
as $(at-ta)q\in (t)^2Q=(t)(s)Q\subseteq (t) R$ (since
$(s)Q\subseteq R$). The claim is obvious in this case.
 Suppose that the second condition holds. Then there are obvious
 $\bA$-module (and $A$-module) isomorphisms:
 $$ (s)Q/(s)R=tQ/tR\simeq {}^{\tau_t}(Q/R)= {}^{\tau_t}M\simeq
 {}^{\btaut}M.$$
Since $\btaut$ is an automorphism of the ring $\bA$, $\pd_{\bA}
({}^{\btaut}M)=\pd_{\bA}(M)=1$. The proof of the claim is
complete. The $\bA$-module $Q/(s)Q$ is projective and
$\pd_{\bA}(M)=1$, hence the $\bA$-module $R/(s)Q$ is projective,
by (\ref{RQM}). Then the short exact sequence (\ref{1RQM}) splits,
$$ R/(s)R\simeq (s)Q/ (s)R\bigoplus R/(s)Q.$$ Therefore,
$\pd_{\bA}(R/(s)R)\geq \pd_{\bA}((s)Q/(s)R)=1$, hence $R$ is {\em
not} a projective $A$-module (otherwise, we would have
$\pd_{\bA}(R/(s)R)=0$, a contradiction), and so $\pd_A(M)>1$ since
$M=Q/R$ and $Q$ is a free $A$-module.  $\Box $

%$\noindent $

%\begin{corollary}\label{d22May10}%\marginpar{d22May10}
%Let $A$ be a ring, $s$ be its regular and normal element such that
%$as-sa\in (s^2)$ for all elements $a\in A$; $\bA= A/(s)$. If $M$
%is  a nonzero $\bA$-module and $\pd_{\bA}(M)=n<\infty$ then
%$\pd_A(M)=n+1$.
%\end{corollary}

%{\it Proof}. Put $t=s$ in Theorem \ref{22May10}.  $\Box $

$\noindent $

Theorem \ref{22May10} is used in the proof of Lemma
\ref{e22May10}.(3) which is the key moment in the proof of Theorem
\ref{AH22May10}.

\begin{lemma}\label{e22May10}%\marginpar{e22May10}
Let $A$ be a ring. Then
\begin{enumerate}
\item For each $\mS_1(A)$-module $M$, $\pd_A(M)\leq
\pd_{\mS_1(A)}(M)$.  \item For each $A$-module $N$,
$\pd_{\mS_1(A)}(\mS_1(A)\t_A N)=\pd_A(N)$.\item If $M$ is a
nonzero $\mS_1(A)$-module such that  $(x-1)M=0$ then
$\pd_{\mS_1(A)}(M) = \pd_A(M)+1$.
\end{enumerate}
\end{lemma}

{\it Proof}. 1. Since $\mS_1(A)$ is a free left $A$-module, any
projective resolution of the $\mS_1(A)$-module $M$ is also a
projective resolution of the $A$-module $M$. Therefore,
$\pd_A(M)\leq \pd_{\mS_1(A)}(M)$.

2. It follows from  the decomposition  $\mS_1(A)
=\bigoplus_{i,j\in \N}x^iy^jA$ that the $A$-module $\mS_1(A)
\t_AN=\bigoplus_{i,j\in \N}Ax^iy^j\t_AN\simeq ({}_AN)^{(\N^2)}$ is
the direct sum of $\N^2$ copies of the $A$-module $N$. Then, by
statement 1, $\pd_{\mS_1(A)}(\mS_1(A)\t_AN)\geq \pd_A(N)$. If
$P\ra N$ is a projective resolution of the $A$-module $N$ then
$\mS_1(A)\t_AP\ra \mS_1(A)\t_AN$ is the projective resolution of
the $\mS_1(A)$-module $\mS_1(A)\t_AN$ since $\mS_1(A)_A$ is free.
Therefore, $\pd_{\mS_1(A)}(\mS_1(A)\t_A N)\leq \pd_A(N)$.

3. Let $n=\pd_A(M)$. If $n=\infty$ then
$\pd_{\mS_1(A)}(M)=\infty$, by statement 1. If $n<\infty$ then we
 will see that the equality  follows from Theorem \ref{22May10} where $s=x-1$ and
$t=y-1$. For, we have to check that the conditions of Theorem
\ref{22May10} hold for the elements $s$ and $t$ (we will see that
both conditions 1 and 2 of Theorem \ref{22May10} hold). By Lemma
\ref{c22May10}, the elements $s$ and $t$ of the ring $\mS_1(A)$
are regular. By Lemma \ref{a22May10}, the element $s$ is left
normal and the element $t$ is right normal with $(s) = (t)$.
Consider the inner derivation $\d =\ad (t)= \ad (y)$ of the ring
$\mS_1(A)$ (recall that $\ad (t)(a)=ta-at$ for all elements $a\in
\mS_n(A)$). Since $\d (A)=0$, $\d (y)=0$ and $\d (x) =
yx-xy=E_{00}$, we see that $\d (\mS_1(A))\subseteq F(A)$. Notice
that $F(A)^2= F(A)$, by (\ref{1SnAs}). By Lemma \ref{a22May10},
$F(A)\subseteq (y-1)$. Then $at-ta \in F(A)=F(A)^2\subseteq
(y-1)^2=(t^2)$ for all elements $a\in \mS_1(A)$. By Lemma
\ref{a22May10}.(2), $\overline{\mS_n(A)}=\mS_n(A)/(t)=A$ and the
ring endomorphism $\btaut$ is the identity automorphism of the
ring $A$. This means that the assumptions of Theorem \ref{22May10}
hold for the elements $s=x-1$ and $t=y-1$, and so
$\pd_{\mS_1(A)}(M) = n+1$, by Theorem \ref{22May10}. $\Box $

$\noindent $

{\bf Proof of Theorem \ref{AH22May10}}. Since $\mS_n(A) =
\mS_1(\mS_{n-1}(A))$, it suffices to prove statements when $n=1$.
In view of the involution $*$ on the algebra $\mS_1(\Z )$ such
that $y^*=x$, and $\mS_1(A)= A\t_\Z\mS_1(\Z )$, it suffices to
show that $\lgldim (\mS_1(A)) = \lgldim (A) +1$.

By Lemma \ref{e22May10}.(2), $\lgldim (\mS_1(A))\geq \lgldim (A)$,
and so the equality is true if $\lgldim (A) = \infty$. Suppose
that $\lgldim (A) < \infty$. By Lemma \ref{e22May10}.(3), $\lgldim
(\mS_1(A)) \geq  \lgldim (A) +1$ since each $A$-module $N$ can be
seen as the $\mS_1(A)$-module which annihilated by the element
$x-1$ since ${}_A\mS_1(A))_A= A\bigoplus (x-1)$ and $(x-1) =
\mS_1(A) (x-1)$. It remains to show that $\pd_{\mS_1(A)}(M)\leq
\lgldim (A) +1$ for all nonzero $\mS_1(A)$-modules $M$.  There is
the short exact sequence of $\mS_1(A)$-modules (by Lemma
\ref{c22May10} and (\ref{S1Z})) %\marginpar{SAM}
\begin{equation}\label{SAM}
0\ra \mS_1(A) \t_A M\stackrel{\cdot (x-1)}{\ra} \mS_1(A)\t_A M\ra
M\ra 0.
\end{equation}
Notice that $\mS_1(A) \t_A M = \mS_1(\Z )\t_\Z A\t_A M\simeq
\mS_1(\Z ) \t_\Z M$, and the map $\cdot (x-1)$ above is equal to
the map
$$ \cdot (x-1)\t {\rm id}_M: \mS_1(\Z ) \t_\Z M\ra \mS_1(\Z )
\t_\Z M. $$ By Lemma \ref{e22May10}, $\pd_{\mS_1(A)}(\mS_1(A)\t_A
M) =\pd_A(M) \leq \lgldim (A)$. By (\ref{SAM}), $\pd_{\mS_1(A)}(M)
\leq \lgldim (A) +1$. The proof of the theorem is complete. $\Box
$

%$\noindent $

\begin{corollary}\label{f22May10}%\marginpar{f22May10}
\begin{enumerate}
\item {\rm \cite{shrekalg}} Let $K$ be a field then $\lgldim
(\mS_n(K))= \rgldim (\mS_n(K))=n$.  \item $\lgldim (\mS_n(\Z ))=
\rgldim (\mS_n(\Z ))=n+1$.
\end{enumerate}
\end{corollary}

\begin{theorem}\label{23May10}%\marginpar{23May10}
Let $A$ be an algebra over a field $K$ such that either $\wdim (A)
= \lgldim (A)$ (eg, $A$ is a left Noetherian algebra) or $\wdim
(A) = \rgldim (A)$ (eg, $A$ is a right Noetherian algebra). Then
$\wdim (\mS_n(A))= \wdim (A) +n$.
\end{theorem}

{\it Proof}. For any two $K$-algebras $R$ and $S$, $\wdim (R\t S)
\geq \wdim (R) +\wdim (S)$  (Theorem 16, \cite{Aus}). Therefore,
$\wdim (\mS_n(A))=\wdim (\mS_n\t A)\geq \wdim (\mS_n)+\wdim (A) =
n+\wdim (A)$ since $\wdim (\mS_n)=n$, Corollary 6.8,
\cite{shrekalg}. The inverse inequality follows from Theorem
\ref{AH22May10} and the fact that $\wdim (R)\leq \min \{ \lgldim
(R), \rgldim (R)\}$ for all rings $R$. Let $d=\lgldim$ (resp. $d=
\rgldim$). If $A$ is a left (resp.  right)  Noetherian algebra
then $d(A)=\wdim (A)$ and
$$\wdim (\mS_n(A)) \leq d(\mS_n(A))=d(A) +n = \wdim (A) +n.$$
The proof of the theorem is complete. $\Box $

%$\noindent $

\begin{corollary}\label{a23May10}%\marginpar{a23May10}
$\wdim (\mS_n(\Z ))=n+1$.
\end{corollary}

In his report the referee asked the following two questions ($K$
is a field):

Q1. {\em Is the algebra $\mS_n(K)$ coherent?}

Q2. {\em What is $K_0(\mS_n(K))$?}

Theorem \ref{31Mar11} gives the answer to the first question. A
module $M$ over a ring $R$ is {\em finitely presented} if there is
an exact sequence of modules $R^m\ra R^n\ra M\ra 0$. A finitely
generated module is a {\em coherent} module if every finitely
generated submodule is finitely presented. A ring $R$  is a {\em
left} (resp. {\em right}) {\em coherent ring} if the module
${}_RR$ (resp. $R_R$) is coherent. {\em A ring $R$ is a left
coherent ring iff, for each element $r\in R$, $\ker_R(\cdot r)$ is
a finitely generated left $R$-module and the intersection of two
finitely generated left ideals is finitely generated}, Proposition
13.3, \cite{Stenstrom-RingQuot}. Each left Noetherian ring is left
coherent but not vice versa.

\begin{theorem}\label{31Mar11}%\marginpar{31Mar11}
Let $K$ be a field and $\mS_n = \mS_n(K)$. Then the algebra
$\mS_n$ is a left coherent algebra iff the algebra $\mS_n$ is a
right coherent algebra iff $n=1$.
\end{theorem}

{\it Proof}. In view of existence of the involution $*$ on the
algebra $\mS_n$
 the first `iff' is obvious since the algebra $\mS_n$ is
 self-dual,
 i.e. it is isomorphic to its opposite algebra
$\mS_n^{op}$ (via $\mS_n\ra \mS_n^{op}$, $a\mapsto  a^*$).
 By Corollary \ref{f22May10}.(1), the algebra $\mS_1$ is a
 hereditary algebra hence coherent.

 If $n\geq 2$ then the algebra
 $\mS_n$ is not left coherent as the left $\mS_n$-module $\CK
 :=\ker_{\mS_n}(\cdot (x_1-x_2))$ is not finitely generated. In
 more detail, $\CK = \CK' \t \mS_{n-2}$ where $\CK'
 :=\ker_{\mS_2}(\cdot (x_1-x_2))$ and $\mS_n = \mS_2\t \mS_{n-2}$.
 It suffices to show that the $\mS_2$-module $\CK'$ is an infinite
 direct sum of nonzero $\mS_2$-modules (then $\CK$ is an infinite
 direct sum of nonzero $\mS_n$-modules, hence $\CK$ is not a
 finitely generated $\mS_n$-module and as a result the algebra
 $\mS_n$ is not left coherent)). So, we can assume that $n=2$.
 Recall that $\mS_2= \mS_1(1)\t \mS_1(2)$  where $\mS_1(k) =
 K\langle x_k, y_k\rangle\simeq \mS_1$ for $k=1,2$; $F(k):=
 \bigoplus_{i,j\in \N} KE_{ij}(k)$ is the only proper ideal of the
 algebra $\mS_1(k)$ and  $\mS_1(k) /F(k)=L_1(k):=K[x_k,
 x_k^{-1}]$. The ideal $\ga_2 := F(1)\t \mS_1(2)+\mS_1(1)\t F(2)$
 of the algebra $\mS_1$ is such that $\mS_2/ \ga_2 \simeq L_1(1)\t
 L_1(2)= K[x_1, x_1^{-1}, x_2, x_2^{-1}]$ is a commutative domain.
  Therefore, $\CK = \ker_{\ga_2} (\cdot (x_1-x_2))$ since the
  element $x_1-x_2$ is a nonzero element of the domain $L_1(1)\t
 L_1(2)$. It follows from the short exact sequence of right
 $\mS_2$-modules
 $$ 0\ra F:= F(1)\t F(2)\ra \ga_2\ra \ga_2/(F(1)\t \mS_1(2)) \simeq
 L_1(1)\t \mS_1(2) \ra 0$$
and the equality $\ker_{L_1(1)\t \mS_1(2)} (\cdot (x_1-x_2))=0$
(use the fact that the algebra $L_1(1)=\bigoplus_{i\in \Z} Kx_1^i$
is $\Z$-graded) that $\CK = \ker_{F} (\cdot (x_1-x_2))$. The ideal
$$F= \bigoplus_{\alpha , \beta \in
\N^2}KE_{\alpha\beta}=\bigoplus_{\beta \in \N^2}\mS_2E_{0\beta}=
\bigoplus_{\alpha \in \N^2}E_{\alpha 0}\mS_2$$ is the direct sum
of nonzero left ideals $\mS_2E_{0\beta}=\bigoplus_{\alpha \in
\N^2} KE_{\alpha \beta}$ and the direct sum of nonzero right
ideals $E_{\alpha 0}\mS_2=\bigoplus_{\beta \in \N^2} KE_{\alpha
\beta}$. Using the identities
$$E_{i,j}(k)x_k=\begin{cases}
E_{i,j-1}(k)& \text{if }j\geq 1,\\
0& \text{if }j=0,\\
\end{cases}
$$
we see that $\ker_{E_{\alpha 0}\mS_2}(\cdot (x_1-x_2)) =
\bigoplus_{s\in \N}Kv_\alpha (s)$ where $v_\alpha (s):=
\sum_{\beta_1+\beta_2=s} E_{\alpha\beta}$ and $\beta = (\beta_1,
\beta_2)\in \N^2$. Notice that $v_\alpha (s) = x^\alpha v_0(s)$
for all $\alpha \in \N^2$ and $y^\beta v_0(s)=0$ for all $\beta
\in \N^2\backslash \{ 0\}$. Therefore, $$\ker_F(\cdot
(x_1-x_2))=\bigoplus_{\alpha \in \N^2} \ker_{E_{\alpha
0}\mS_2}(\cdot (x_1-x_2))=\bigoplus_{\alpha\in \N^2, s\in
\N}Kv_\alpha (s) = \bigoplus_{s\in \N} \mS_2v_0(s)$$ is an
infinite direct sum of nonzero left ideals of the algebra $\mS_2$,
as required.  $\Box $

$\noindent $

$${\bf Acknowledgements}$$

The author would like to thank the referee  for the interesting
comments and challenging questions.

\small{

Department of Pure Mathematics

University of Sheffield

Hicks Building

Sheffield S3 7RH

UK

email: v.bavula@sheffield.ac.uk

}

\end{document}